\documentclass[12pt]{amsart}
\usepackage{amsmath,amsfonts,amsthm,amsopn}
\usepackage{graphicx}
\usepackage{epsfig,verbatim}
\usepackage{xspace}
\usepackage{mathrsfs}
\newtheorem{Def}{Definition}
\newtheorem{Cor}{Corollary}
\newtheorem{Th}{Theorem}
\newtheorem{Pro}{Proposition}
\newtheorem{Lem}{Lemma}
\newtheorem{ex}{Example}
\theoremstyle{remark}
\newtheorem{rem}{Remark}
\def\Wsp{{\boldsymbol W}} 
\newcommand{\Lz}{\mathbf{L}^2 (\mathbb{R}^d )}
\newcommand{\Rzd}{\mathbb{R}^{2d}}
\newcommand{\Gmgr}{\mathcal{G}^Q_{m}}
\newcommand{\Id}{\operatorname{Id}}
\begin{document}

\title[Quilted Frames]{Quilted Gabor frames - a new concept for adaptive time-frequency representation}

\author{Monika D\"orfler}
\address{Institut f\"ur Mathematik, Universit\"at Wien, Alserbachstrasse 23
  A-1090 Wien,   Austria }
\email{monika.doerfler@univie.ac.at} %\subjclass{}
\thanks{M.\ D.\ was supported by the FWF grant T~384-N13}
\subjclass[2000]{ 42A65,42C15,42C40}
\date{\today}
\keywords{Time-frequency analysis \, adaptive representation \,  uncertainty principle \,  frame bounds \,  frame algorithm} 
\begin{abstract}
Certain signal classes such as audio signals call for signal representations with the ability to adapt to the signal's properties.
In this article we introduce the new concept of quilted frames, which aim at adaptivity in time-frequency representations. As opposed to Gabor or wavelet frames, this new class of frames allows for the adaptation of the signal analysis to the local requirements of signals under consideration.  Quilted frames are constructed  directly in the time-frequency domain in a signal-adaptive manner. Validity of the frame property guarantees the possibility to reconstruct the original signal. The frame property is shown  for specific situations and the Bessel property is proved for the general setting. 
Strategies for reconstruction from coefficients obtained with quilted Gabor frames and numerical simulations are provided as well.  
\end{abstract} 

\maketitle

\section{Introduction}
\begin{center}\begin{scriptsize}
QUILT (verb): (a) to fill, pad, or line like a quilt\\
(b) to stitch (designs) through layers of cloth \\(c) to fasten between two pieces of material
\end{scriptsize}\end{center}
Natural  signals usually  comprise components of various different characteristics and their analysis requires judicious choice of processing tools.  For audio signals time-frequency dictionaries have proved to be an adequate option. Since orthonormal bases cannot provide good time-frequency resolution,~\cite{gr01}, time-frequency  analysis naturally leads to the use of {\it frames}.
Most  classes of frames commonly used in applications, be it wavelet or Gabor
frames, feature a resolution following a fixed rule over the whole
time-frequency or time-scale plane, respectively. The concept of {\it quilted frames}, as introduced in this contribution,  gives up this uniformity and allows for
different resolutions in assigned areas of the time-frequency
plane.\\
 The primary motivation for introducing this new class of frames stems from the processing of audio
and in particular music signals, where the trade-off between time- and frequency resolution has a strong impact on the results of
analysis and synthesis, see \cite{ro98,do01,WDG01,prruwo09,badoja09}.  The well-known uncertainty principle makes the choice of just one 
analysis window a difficult task:  different resolutions might be
favorable in order to achieve sparse and precise representations for the various signal components. For example,  percussive
elements  require short analysis windows and high sampling rate in
time, whereas sustained sinusoidal components are better
represented with wide windows and a longer FFT, thus more sampling
points in frequency.

Several approaches have been suggested to deal with the trade-off
in time-frequency resolution. The notion of  \emph{multi-window}
Gabor expansions, introduced by Zibulski and Zeevi,~\cite{zezi97},	
uses a finite number of windows of different shape  in order to
obtain a richer dictionary with the ability to better represent
certain characteristics in a given signal class. Another approach
is the usage of several \emph{bases} in order to best describe the
components of a signal with a priori known characteristics, see
\cite{DauTor}. All these approaches, however, stick to a uniform
resolution guided by the action of a certain group via a unitary
representation. 
For quilted Gabor frames we give up this restriction and
introduce systems  constructed from globally defined
frames by restricting these to
certain, possibly compact, regions in the time-frequency or
time-scale  plane. The  idea of realizing tilings of the time-frequency plane has been
suggested in \cite{Kovac99} and \cite{Vul03}, however, these
authors stick to the construction of orthogonal bases. In this case, every basis function corresponds to a particular tile.   We
will achieve a wider range of possible partitions, windows and sampling
schemes by allowing for redundancy.  Thus we aim at 
designing systems that can optimally adapt to a class of signals
considered. As a particular example of quilted frames, the notion of  reduced multi-Gabor frames
was first introduced in~\cite{do02} and  recently exploited in~\cite{badoja09}.
Note that this model allows, for example, a transform yielding constant-Q spectral resolution, which is invertible, as opposed to the original construction~\cite{BR91}.
Reduced multi-Gabor frames were successfully applied to
the task of denoising corrupted audio signals, see~\cite{WDG01}. The
processing of sound signals also yields a motivation for the
next step in generalizing the idea to \emph{quilted frames}, which
allow arbitrary tilings of the time-frequency plane, see~\cite{JaTo07}. \\
Quilted frames also bear theoretical interested in themselves and should be compared to constructions such as fusion frames~\cite{cakuli08}
and the frames proposed in~\cite{alcamo04-1}. In fact, the construction of quilted frames provides constructive examples for the models presented in these contributions.\\
For the mathematical description of quilted frames, we  start
from principles  of Gabor analysis~\cite{FeSt98}.
%particularly suited to the processing of audio signals, since it employs
%frequency-modulated  building blocks.
The idea for the construction of  quilted Gabor frames is inspired by the early work of Feichtinger and Gr\"obner on decomposition methods~\cite{fegr85,fe87} and recent results on time-frequency partitions for the characterization  of function spaces~\cite{dofegr06,dogr09}:\\
Assume that a covering
%\footnote{The central condition for admissibility is  that $\bigcup_{r \in \mathcal{I}}\Omega _{r} = \mathbb{R}^d$ and
% that the number of overlapping
%areas in the covering  is bounded above. A precise definition is given below.} 
$(\Omega _{r})_{r \in \mathcal{I}}$
of the phase space $\mathbb{R}^{2d}$ is given. To each member of the covering a frame from a family of 
 Gabor frames is assigned, hence, the  new system locally resembles the original frames.
The resulting  global system  will be called a
 \emph{quilted Gabor system}. We conjecture that these systems may be shown to constitute frames under certain, rather general conditions. 
 In this paper we will show the frame property  in two special cases and proof the existence of an upper frame bound for a  general setting.

The rest of this paper is organized as follows. Section~2 provides
notation and gives an overview over basic results in Gabor
analysis. Section~3 introduces the general concept of quilted Gabor frames. In Section~4, the existence of an upper frame bound (Bessel property) for general quilted frames is proved. In Section~5 and Section~6, a lower frame bound is constructed for two particular cases, namely, the partition of the time-frequency plane in stripes and the replacement of frame elements in a compact region of the coefficient domain. Finally, Section~7 presents numerical examples for these cases.
%%%%%%%%%%%%%%%%%%%%%%%%%%%%%%%%%%%%%%%%%%%%%%%%%%%%%%%%%%%%%%%%%%%%%%%%%%%%%%%%
\section{Notation and some basic facts from Gabor theory}
%In this  section we fix some notation and collect 
%important facts from  Gabor theory used in the sequel.

We use  the  normalization 
$\hat{f}(\omega ) = \int_{\mathbb{R}^d} f(t) e^{-2\pi i \omega t}dt$ of the \textbf{Fourier transform} on  $\mathbf{L}^2
(\mathbb{R}^d)$.
 $M_{\omega }$ and $T_x$ denote frequency-shift by
$\omega $ and time-shift by $x$, respectively, of a function $g$,
combined to the \textbf{time-frequency shift operators} $\pi
(\lambda)= M_{\omega } T_x g(t) = e^{2\pi i t\omega }g(t-x)$
for $\lambda = (x,\omega )\in \mathbb{R}^{2d}$.

The \textbf{Short-time Fourier transform}
  (STFT)  of a function $f\in\mathbf{L}^2
(\mathbb{R}^d) $ with respect to a window function
$g\in\mathbf{L}^2 (\mathbb{R}^d) $ is defined  as
\begin{equation}
  \label{Def:STFT}
  \mathcal{V}_g f(\lambda)  = \int _{ \mathbb{R}^d } f(t) \bar{g} (t-x) e^{-2\pi
  i \omega \cdot t} \, dt =  \langle f,   \pi (\lambda ) g\rangle \,
  .
\end{equation}
A {\it lattice} $\Lambda \subset \mathbb{R}^{2d} $ is a discrete
subgroup of $\mathbb{R}^{2d}$ of the form $\Lambda =
A\mathbb{Z}^{2d}$, where $A$ is an invertible $2d\times 2d$-matrix
over $\mathbb{R}$. The special case
$
\Lambda = \alpha\mathbb{Z}^d\times \beta\mathbb{Z}^d ,
$ where $\alpha , \beta > 0$ are the lattice
constants,  is called a separable or product lattice. 
%Unless
%otherwise stated,  we always  assume that the lattice $\Lambda$ in
%use is a product lattice. %The {\em adjoint} lattice
%$\Lambda^{\circ} $ of a product lattice is given by
%\begin{equation}\label{def:adjL}
%\Lambda^{\circ}  = \frac{1}{\beta}\mathbb{Z}^d\times
%\frac{1}{\alpha}\mathbb{Z}^d
%\end{equation}

%\textbf{Frames, Gabor frames, the frame operator and the dual
%frame} \\

A family of functions
$(g_k)_{k \in\mathbb{Z}}$ in $ (\mathbb{R}^d )$ is called a frame,
if there exist lower and upper frame bounds $A,B >0$, so that
\begin{equation}\label{framecond}
A\| f\|^2 \leq \sum_{k \in\mathbb{Z}}|\langle f,g_k\rangle| ^2\leq B\|
f\| ^2\ \mbox{ for all } f\in \mathbf{L}^2 (\mathbb{R}^d)
\end{equation}
Assumption \eqref{framecond} can be understood as an ``approximate
Plancherel formula''. It guarantees that  any
signal $f\in \mathbf{L}^2 (\mathbb{R}^d)$ can be represented as  infinite series
with square integrable coefficients using the elements $g_k$. The existence of the upper bound $B$ is called {\it Bessel property} of the sequence $(g_k)_{k \in\mathbb{Z}}$.
The {\it frame operator} $S$, defined as
\[Sf = \sum_{k\in\mathbb{Z}} \langle f,g_k\rangle g_k \]
allows the calculation of the \textit{canonical dual frame} $(\gamma_k)_{k \in\mathbb{Z}}= (S^{-1}g_k )_{k \in\mathbb{Z}}$, which
guarantees minimal-norm coefficients in the expansion
\begin{equation}\label{fexp}
f = \sum_{k\in\mathbb{Z}} \langle f, \gamma_k\rangle g_k = \sum_{k\in\mathbb{Z}} \langle f,
g_k\rangle \gamma_k .
\end{equation}
If $A=B$, the frame is called {\it tight} and $f = \frac{1}{A}\sum_k \langle f, g_k\rangle g_k$.
We refer the interested reader  to  Christensen's book \cite{Chr03} for more details on general frames.\\
In  the special case  of Gabor or
Weyl-Heisenberg frames,  the frame elements are generated by time-frequency shifts of a basic atom or window $g$ along a  lattice $\Lambda$: 
\[g_\lambda =\pi (\lambda )g.\] In this case we write $S = S_g$, and $S_g
=T_g^{\ast}T_g $, where  $T_g :
f\mapsto [\langle f , g_{\lambda}\rangle ]_{\lambda}$ is  the {\it analysis operator}  mapping the function
$f\in\mathbf{L}^2 (\mathbb{R}^d ) $ to its coefficients $c(f)(\lambda ) = T_g f (\lambda)$. 
These coefficients correspond to the samples of the STFT on $\Lambda$. Its adjoint $T_g^{\ast}:[c_{\lambda} ]_{\lambda\in\Lambda}\mapsto
\sum_{\lambda\in\Lambda}c_{\lambda}\pi(\lambda )g$  is  the
{\it synthesis operator}.
For the Gabor frame generated by time-frequency  shifts of
the window $g$ along the lattice $\Lambda$ we write
$\mathcal{G}(g,\Lambda )$.

We next introduce the concept of {\it
partitions of unity }.  A family $(\psi _{r })_{r \in
\mathcal{I}}$ of non-negative functions with $ \sum_{r}  \psi_r
(x) \equiv 1$ is called \textbf{bounded admissible partition of
unity} (BAPU) subordinate to $ (B_{R_r}(x_r))_{r \in \mathcal{I}}
$, if the support $\Omega _{r}$ of $(\psi _{r})$ is contained in $
B_{R_r}(x_r)$ for $r \in \mathcal{I}$, and $ (  B_{R_r}(x_r) )_{r
\in \mathcal{I}} $ is
 an {\it admissible covering} in the sense of~\cite{FeiGroeb85}
i.e.,\xspace $\bigcup_{r \in \mathcal{I}}B_{R_r} = \mathbb{R}^d$ and
 the number of overlapping
   $B_{R_r}(x_r)$ is bounded above (admissibility condition). In other words, with 
   \[r^{\ast} := \{s : s\in \mathcal{I}, B_{R_r}(x_r)\cap B_{R_s}(x_s )
\neq 0\},\] for all $r \in \mathcal{I}$ there exists $n_0\in\mathbb{N}$, called height of the BAPU,  such that 
$|r ^{\ast}|\leq n_0$.

 For technical reasons, which do not eliminate any interesting
 example, we assume even more:
 for all $\rho<\infty$ the family
$ (  B_{R_r + \rho}(x_r) )_{r \in \mathcal{I}} $
  should be an admissible covering of $\mathbb{R}^d$.
More precisely, we assume throughout this paper that for each
$\rho > 0 $
   there exists $n_0 = n_0(\rho) \in\mathbb{N}$ such that
   the number of overlapping balls constituting the covering
   is uniformly controlled:
$|r ^{\ast}|\leq n_0$ for all $r \in \mathcal{I}$, where \[r
^{\ast} := \{s : s\in \mathcal{I}, B_{R_r+ \rho}(x_r)\cap
B_{R_s}(x_s + \rho) \neq 0\}.\] Obviously such coverings are of
uniform height.\\
Using the concept of BAPUs, we now turn to 
\textbf{Wiener amalgam spaces}, introduced by H.~Feichtinger
in 1980 (see~\cite{fe92-3} for an accessible publication).
The definition of Wiener amalgam spaces aims at decoupling  local and global properties
of $\mathbf{L}^p$-spaces. 
Let a BAPU  $(\psi _{r })_{r \in
\mathcal{I}}$ for $\mathbb{R}^d$ be given. 
The Wiener amalgam space $\Wsp(\mathbf{L}^p,\ell^q )$ is defined as
follows:
\[\Wsp(\mathbf{L}^p ,\ell^q )(\mathbb{R}^d) = \big\{ f \in L^p_{loc} : \| f\|_{\Wsp(\mathbf{L}^p ,\ell^q )} =
\big(\sum_{r \in
\mathcal{I}} \| f\cdot
 \psi_r\|_p^q\big)^{\frac{1}{q}}<\infty\big\}.\]
We will denote by $\Wsp(C^0, \ell^p )(\mathbb{R}^d)\subseteq \Wsp(\mathbf{L}^\infty ,\ell^p )(\mathbb{R}^d)$ the subspace of continuous, locally bounded functions in  $\Wsp(\mathbf{L}^\infty ,\ell^p )(\mathbb{R}^d)$. A
 comprehensive review of (weighted) Wiener amalgam spaces can be found	
in~\cite{Heilamal}. We note that in their most general form they
are described as $\Wsp(B,C)$, with local component $B$ and global
component $C$. Let us recall some properties  which will be needed
later on:
\begin{itemize}
\item If $B_1\subseteq B_2$, $C_1\subseteq C_2 $, then $\Wsp( B_1 ,
C_1 ) \subseteq \Wsp( B_2 , C_2 ) $.
%\item $\mathbf{L}^2  \bigcup \mathbf{L}^\infty\subseteq \Wsp(\mathbf{L}^\infty , \ell^2 )$.
\item If $B_1\ast B_2\subseteq B_3$, $C_1\ast C_2\subseteq C_3$,
then $\Wsp(B_1 ,C_1 )\ast \Wsp(B_2 , C_2 )\subseteq \Wsp(B_3 , C_3 )$.
\end{itemize}
A  particularly important Banach space  in time-frequency analysis is the Wiener Amalgam space $\Wsp(\mathcal{F}L^1,\ell^1)$. This space, also known under the  name Feichtinger's algebra, is better known as the modulation space   $M^{1,1}_m $, 
   with
constant weight $m \equiv 1$. It is often denoted by $S_0$ in the literature and we will adopt this name in the present work. For convenience, we also recall the definition of $S_0$ via the short-time Fourier transform. 
\begin{Def}[$S_0$] Let
  $g_0$ be the Gauss-function $g_0 = e^{-\pi\| x\|^2}$.
The  space $S_0 (\mathbb{R}^d)$ is given by
\[S_0( \mathbb{R}^d)
=\{f\in\Lz :\|f\|_{S_0}=\|\mathcal{V}_{g_0}f \|_{\mathbf{L}^1
(\Rzd )}<\infty\}.\]
\end{Def} An in-depth investigation of $S_0(\mathbb{R}^d)$
and its outstanding role in time-frequency  analysis
 can be found in~\cite{FZ98a}. Note that  $S_0(\mathbb{R}^d)$ is densely
embedded in $\mathbf{L}^2 (\mathbb{R}^d) $, with $\| f\|_2 \leq \| g_0\|_2^{-1}\|
f\|_{S_0}$. Its dual space $S_0'(\mathbb{R}^d)$, the space of all linear,
continuous functionals on $S_0(\mathbb{R}^d)$, contains $\mathbf{L}^2(\mathbb{R}^d) $ and is a very
convenient space of (tempered) distributions.
 Moreover, in the definition of  $S_0(\mathbb{R}^d) $, $g_0$ can be replaced by
any $g\in S_0(\mathbb{R}^d)$, see~\cite[Theorem 11.3.7]{gr01} and different
functions  $g \in S_0(\mathbb{R}^d) \setminus \{0\} $ define equivalent norms on
$S_0(\mathbb{R}^d)$.\\
One of the results of major importance in Gabor analysis states
that for $g\in S_0(\mathbb{R}^d)$ the analysis mapping $T_g$ is bounded from
$\mathbf{L}^2(\mathbb{R}^d) $ to $\ell^2$ for \emph{any} lattice $\Lambda$, and
$T^{\ast}_g$ is then bounded by duality, see \cite[Section
3.3]{FZ98a} for details. This will be of crucial significance in
our arguments.

%%%%%%%%%%%%%%%%%%%%%%%%%%%%%%%%%%%%%%%%%%%%%%%%%%%%%%%%%%%%
\section{Quilted Gabor frames: the general concept}
\begin{figure}[tb]
\centerline{\includegraphics[scale = .8]{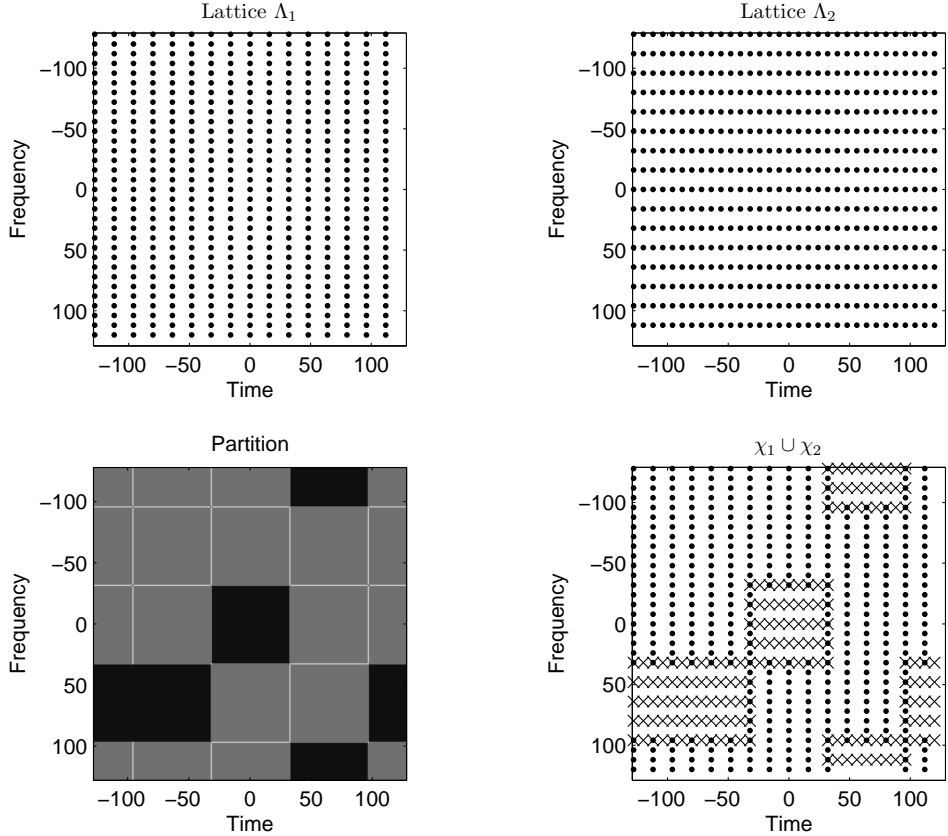}}
\caption{{\it Partition in time-frequency and resulting quilted lattice}}\label{Fig1}
\end{figure}

For  the construction of quilted frames, we  start from a collection of (Gabor) frames. Usually, these frames will feature various different qualities, e.g. varying resolution quality for time and frequency. Then, a partition in time-frequency is set up according to some application-dependant  criterion and a particular frame is assigned to each member of the partition. For example, in~\cite{JaTo07}, the selection of the local frames is based on time-frequency sparsity criteria. 
Figure~\ref{Fig1} gives an illustration of the basic idea, for a partition assigning one out of two  different Gabor frames to each of the tiles of size $64\times 64$, where the signal length is $L = 256$ and the number of  tiles thus	$16$. The upper displays show the lattices corresponding to the two Gabor frames, the last display shows the ``quilted" lattice $\chi_1\cup\chi_2$ resulting from concatenation. Let us emphasize at this point, that at sampling points marked with different symbols, different windows are also used. 
%We will come back to this example in Section~\ref{Se:ExSim}.\\
Note that, conceptually, irregular tilings may be used just as well. However, for practical as well as theoretical reasons, tilings with some kind of structure are more 	beneficial.\\
It is important to point out that the partition in different domains corresponding to various different frames actually happens \emph{in the time-frequency domain}. This implies that a priori we have no knowledge about the properties of the local families, as opposed to the concept of fusion frames, as discussed in~\cite{cakuli08,caku04}. In particular, we are not necessarily dealing with closed subspaces which may be transformed into each other as in the approach introduced in~\cite{fo04}. 

We now give a precise definition for quilted Gabor frames. 
\begin{Def}[Quilted Gabor frames]Let Gabor frames 
$ \mathcal{G}(g^j, \Lambda ^j )_{j\in\mathcal{J}}$  for $\mathbf{L}^2 (\mathbb{R}^d )$ and an admissible covering
$  \Omega_r \subseteq (  B_{R_r}(x_r) )_{r \in \mathcal{I}} $ of $\mathbb{R}^{2d}$ be given. Define the local index sets $\mathcal{X}^r = \Omega_r\cap \Lambda^{m(r)}$, where $m: \mathcal{I}\mapsto \mathcal{J}$ is a mapping assigning a frame from the given Gabor frames to each member of the covering. 
Then the set 
\begin{equation}\label{Def:QF}
\bigcup_{r\in\mathcal{I}} \mathcal{G}(g^{m(r)}, \mathcal{X}^r )	
\end{equation}
is called a quilted Gabor frame for $\mathbf{L}^2 (\mathbb{R}^d )$, if there exist constants $0< A, B < \infty$, such that 
 \begin{equation}\label{Def:QF1}
A\| f\|_2^2\leq\sum_{r\in\mathcal{I}}\sum_{\lambda \in \mathcal{X}^r} |\langle f, \pi (\lambda ) g^{m(r)}\rangle |^2\leq B \| f\|_2^2
\end{equation}
holds for all $f\in\mathbf{L}^2 (\mathbb{R}^d )$. 
 \end{Def}
The general setting of quilted frames  includes, of course, various special cases. We first give some trivial examples which may however be relevant in applications.
\begin{ex}\label{Ex00}
For a given Gabor frame, we may choose additional sampling points in any selected region. This may be helpful, if in some applications, finer resolution is only desirable in certain parts of the time-frequency domain. 
Formally, this may be rephrased as follows. We are given Gabor frames  $\mathcal{G}(g, \Lambda^j )_{j\in\mathbb{N}}$  for $\mathbf{L}^2 (\mathbb{R}^d )$ with $\Lambda^0 \subseteq \Lambda^j$ for $j\in\mathbb{N}$ and $A_0$ the lower frame bound for $j= 0$. Then, for an admissible covering, the local index sets are defined by $\mathcal{X}^r = \Omega_r\cap \Lambda^{m(r)}$, where $m: \mathcal{I}\mapsto \mathbb{N}$ is the mapping selecting the local systems. It is then trivial to see, that the resulting quilted Gabor frame has a lower frame bound $A_0$. The existence of an upper frame bound is covered by Theorem~\ref{mainth}.
\end{ex}

\begin{ex}\label{Ex01}For a given  multi-window Gabor frame $\mathcal{G}(\{g_1,\ldots , g_n\}, \Lambda^0 )$, additional sampling points for selected windows may be added in certain parts of the time-frequency domain. For a formal description, assume that an admissible covering $\Omega_r$, $r\in\mathcal{I}$, is given and let  $\Lambda^0 \subseteq \Lambda^j$ for $j = 1,\ldots, N$ as in the previous example. The mapping $m:\mathcal{I}\rightarrow \{ 0,\ldots, N\}^n$ is given by $m(r )(k) = 0,\, k = 1,\ldots , n$ whenever the original lattice is maintained for all windows in the support of $\Omega_r$ and by $m(r)(k) = j_0, j_0 \in \{ 1,\ldots, N\}$, if denser sampling correspondong to $\Lambda^{j_0}$ is desired in $\Omega_r$ for the window $g_k$.  
\end{ex}

In the next section we will prove the Bessel property of quilted systems obtained in a rather general situation, allowing for a finite overlap between the local patches. In the construction of lower frame bounds, a certain overlap between adjoint patches is often necessary. The  two subsequent sections then describe two situations, in which a lower frame bound for the resulting quilted Gabor frame can be constructed explicitly.

%%+++++++++++++++++++++++++++++++++++++++++++++++++++++++++++++++++++++++++++++++++++++++++++++++++++++++++++++++++++++
\section{The Bessel condition in the general case}\label{Se:BesselCond}

We prove  the existence of an upper frame bound for  quilted frames as defined in \eqref{Def:QF}.  Note that the Bessel property alone allows for interesting conclusions about operators associated with the respective sequence, compare~\cite{ba07}.
We will deduce the Bessel property of quilted Gabor frames from a general statement on relatively separated sampling sets. This result generalizes a result given in~\cite{lawa03} on the Bessel property of irregular time-frequency shifts of a single atom.  We prove that an arbitrary function from a set of window functions satisfying a common decay condition may be chosen for every sampling point in a relatively separated sampling set to obtain  a Bessel sequence.

 We assume that
different given Gabor systems are to be used in  compact sets $\Omega_r$
corresponding to the members of an admissible covering of $\mathbb{R}^{2d}$.  Under the
assumption that the windows under consideration satisfy a common decay
condition in time-frequency and that the set of  lattices is compact,
we claim that an upper frame bound, or Bessel bound, can be found.  As before, $g_0$ denotes the Gaussian window. 
 \begin{Th}\label{mainth} 
Let  Gabor frames $ \mathcal{G}(g^j, \Lambda ^j )_{j\in\mathcal{J}}$  for $\mathbf{L}^2 (\mathbb{R}^d )$
 and an admissible covering
$  \Omega_r \subseteq (  B_{R_r}(x_r) )_{r \in \mathcal{I}} $
 of the signal domain be given.  Assume further that
\begin{itemize}
\item[(i)] $g^j\in H_{s,C}$ for all $j$, where
\[H_{s,C} = \{g\in \mathbf{L}^2(\mathbb{R}^d)  :|\mathcal{V}_{g_0} g| (z)\leq C  (1+\|z\|^2)^{-\frac{s}{2}}\},\ s>2d,\ C>0 .\]
\item[(ii)] the lattice constants $\alpha ^j,\beta ^j$  are chosen from a compact
set in $\mathbb{R}^+\times\mathbb{R}^+$, \\ i.e. $\alpha^j
\subseteq [\alpha_0, \alpha_1] \subset (0,\infty) $ and $\beta^j
\subseteq [\beta_0, \beta_1] \subset (0,\infty). $
\item[(iii)]The regions assigned to the different Gabor
systems correspond to an admissible covering $\Omega_r, r \in\mathcal{I}$ with 
 $\operatorname{supp} (\psi _r )\subseteq  \Omega_r \subseteq B_{R_r}(x_r)  $
 for $r\in\mathcal{I}$. 
 \end{itemize}
 Let $m: \mathcal{I}\mapsto \mathcal{J}$ be a mapping assigning a frame from $ \mathcal{G}(g^j, \Lambda ^j )_{j\in\mathcal{J}}$ to each member of the covering. 
Then for any $\delta<\infty$, the overall family given by 
\begin{equation}\label{mixsyst1}
\Gmgr = \bigcup_{r\in \mathcal{I}}\{ \pi (\lambda )g^{m(r)}:  \lambda\in
\mathcal{X}^r\subset\Lambda^r, \ \mathcal{X}^r= \Lambda^{m(r)}\cap B_{R_r +\delta} \}
\end{equation} possesses an upper frame bound, ~i.e., is a Bessel sequence for $\Lz$.
\end{Th}
Note that the theorem states that \emph{in particular} the local systems given by $\mathcal{X}^r = \Omega_r\cap \Lambda^{m(r)}$ for all $r$ lead to a Bessel sequence. More generally, however, the local patches can uniformly be enlarged by a radius $\delta$.\\
We first prove a general statement on
sampling of functions in certain Wiener amalgam spaces 
over   relatively separate sampling sets. 
\begin{Def}[Relatively separated sets]
A set $\mathcal{X} = \{z_i = (x_i, \xi _i ), i\in\mathcal{I}\}$ in $\mathbb{R}
^{2d}$ is called (uniformly) $\gamma$-separated, if
$\inf_{j,k\in\mathcal{I}, i\neq j}|z_j -z_k |=\gamma>0$. A relatively separated set
is  a \emph{finite} union of 
separated sets. We call $\mathcal{X}$ $(\gamma,R)$-relative separated if the number of separated sets is $R$.
\end{Def}
\begin{rem}
It is  easy to show that the concept of relative separation does
{\it not} depend on the specific values of $\gamma$ and $R$. In
other words, any  $(\gamma,R)$-relative separated set is also a
finite union of $\eta$-separated subsets. Of course one has to
allow to compensate the smallness of $\eta$ by a larger number $R'
= R'(\eta)$.\\
 There is an equivalent point of view. A sequence is
relatively separated in $R^k$ if and only if for some fixed  $s >
0$ the family $(B_s(x_k))_{k \geq 1}$ covers each point $x$ in
$R^k$ at most $h = h(s)$ times, uniformly with respect to $x$.
\end{rem}

\begin{Lem}\label{amalest1}
Let $1\leq p<\infty$.
Given a pair $(\gamma, R)$ there exists a constant
$C = C(\gamma,R)$ such that for all $(\gamma,R)$-separated
sets $\mathcal{X}$ and all
 functions $F \in \Wsp(C^0, \ell^p )(\mathbb{R} ^{2d})$
 one has $F|_{\mathcal{X}}$ is in $\ell ^p (\mathcal{X} )$ and
\begin{equation}\label{amalnormrelsep}
 \|F |_{\mathcal{X}} \|_{  \ell ^p (\mathcal{X} ) }  =
\left ( \sum_{x_i\in\mathcal{X}} |F(x_i )|^p \right ) ^{1/p}
\leq C  \|F \| \
_{\Wsp(C^0, \ell^p )}.
\end{equation}
\end{Lem}
\underline{Proof:}\\ Recall that $\|F\|_{\Wsp(C^0, \ell^p )} =\|
a_{kn}\|_{\ell^p}$, where
\[a_{kn} = \mathrm{esssup}_{(x,\xi )\in Q}|F(x+k, \xi+n)|
=\|F\cdot T_{(k,n)}\chi_{Q}\|_{\infty}.\] By assumption,
$\mathcal{X}$ is the finite union of uniformly separated sets $\mathcal{X}^r, r = 1,\ldots ,
R$ and there exists $\gamma$, such that $\min|z_j -z_k |\geq\gamma
>0$ for any pair $z_j, z_k$ in $\mathcal{X}^r$, with $j \neq k$. Hence,
for all $ r = 1,\ldots , R$, we have at most $  (1 +
\frac{1}{\gamma})^{2d}  $ points $z_i\in\mathcal{X}^r$ in the box $(k,n)+Q$,
$(k,n)\in\mathbb{Z}^d\times\mathbb{Z}^d$, such that the number of
$z_i\in\mathcal{X}$ in  this box is bounded above by $R\cdot(  1 +
1/{\gamma})$.  Clearly
\[|F(x,\xi)|\leq \|F\cdot T_{(k,n)}\chi_{Q}\|_{\infty} \mbox{ for all } (x,\xi)\in (k,n) + Q.\]
Altogether, this yields:
\begin{align}
\sum_{z_i\in\mathcal{X}} |F(z_i )|^p \ \leq  \ & \ R\cdot (1 +
\frac{1}{\gamma})^{2d} \sum_{(k,n)\in\mathbb{Z}^{2d}}\|F\cdot
T_{(k,n)}\chi_{}\|_{\infty}^p\\
= \ & R \cdot (1 +
\frac{1}{\gamma})^{2d}  \ \|F\|^p_{\Wsp(C^0, \ell^p
)}\notag
\end{align}
\hfill$\square$
\begin{rem}
Analogous statements  are true in more general situations,
especially for any weighted sequence space $\ell^p_m$, see
\cite[Lemma 3.8]{fg89jfa} for example.\\
The condition, that $F$ is continuous (locally in $C^0$), guarantees, that sampling is well-defined. Of course,  weaker conditions, for example, semi-continuity, are sufficient.  
\end{rem}
 The upper
frame bound estimate will follow from  a point wise estimate over the family of short-time
Fourier transforms
\begin{equation}
\label{point-max1}
 F(\lambda) =  \sup_{j \in \mathcal{J}} |V_{g^j} f(\lambda)|, \quad \lambda \in
 \mathbb{R} ^{2d}
\end{equation}
We will make use of
following lemma.
\begin{Lem}\label{Le:supest} Assume that
 $g\in M \subseteq H_{s,C}$ for $ s > 2d$, $C > 0$. Then there exists some constant
$C > 0$ such that for all  $f \in \mathbf{L}^2(\mathbb{R}^d) $ one has
the following uniform estimate of $V_{g}f$ in $\Wsp(C^0,\ell^2)$:
\[  \|\sup_{g\in M} |V_{g}f| \|_{\Wsp(C^0, \ell^2)   } \quad \leq
    \quad  C_1  \| f \|_2, \quad \mbox{ for all }  f \in L^2.   \]
\end{Lem}
\underline{Proof:}\\
The crucial step is to invoke the
convolution relation between different short-time
Fourier transforms (\cite[Lemma 11.3.3]{gr01} ).
For convenience in the application below let
us denote the generic elements from $M$ by
$g^j$ (instead of $g_0$ in 11.3.3), and make the choice
 $\gamma = g = g_0$, the normalized Gauss function $g_0$.
Then obviously  $\langle \gamma, g\rangle = \|g_0\|_2^2 = 1$
and we have the following estimate
\[|\mathcal{V}_{g^j} f (\lambda )|\leq
\left[|\mathcal{V} _{g_0} f|\ast |\mathcal{V}_{g^j} g_0 |\right](\lambda )
 \]
 Since $ | \mathcal{V}_{g^j} g_0  (\lambda)| = |\mathcal{V}_{g_0}{g^j}  (-\lambda)| \leq
   w_s(\lambda)$, we may
take the pointwise supremum over $g^j \in M$ on the left side
and arrive at
\[ \sup_{g^j \in M} |\mathcal{V}_{g^j} f (\lambda )|\leq
(|\mathcal{V} _{g_0} f| \ast C w_s)  (\lambda)  \]
 Of course  $s > 2d$  implies that
    $ w_s \in \Wsp(C^0,\ell^1)(\mathbb{R}^{2d}) $.
   Using the general fact  (\cite[Cor. 3.2.2]{gr01} ) that
  $ \| V_{g_0}f\|_2  \leq  C_2 \|f\|_2$ for any $f \in \mathbf{L}^2(\mathbb{R}^d) $
  and applying the convolution relation $\mathbf{L}^2 \ast \Wsp(C^0,\ell^1)
   \subseteq \Wsp(L^1,\ell^2) \ast \Wsp(C^0,\ell^1) \subseteq
     \Wsp(C^0,\ell^2)$, together with the appropriate estimates,
     we arrive at the desired estimate.\hfill$\square$\\
\begin{rem}
It is worthwhile to note that the above result is not just a
simple compactness argument. As a matter of fact it is not
difficult to construct compact sets  $ M \subset S_0(\mathbb{R}^d)$ for which
the above result is not valid. One may, for example, just take a
null sequence of the form $ (c_n T_{x_n} g)_{n \geq 1} $ for some
$g \in S_0(\mathbb{R}^d)$, and with $(c_n) \in C_0$ but $(c_n) \notin \ell^2$.
\end{rem}
The next theorem states that for a relatively separated sampling set of time-frequency shifts we can construct a Bessel sequence by associating to each sampling point an element from $H_{s,C}$. As before, we 
let $M\subseteq H_{s,C}$ be a set of window functions indexed by $\mathcal{J}$. 

\begin{Th}\label{Th:uppbBessel}
Let $\mathcal{X}=\{z_i = (x_i, \xi _i ), i\in\mathcal{I}\}$ in $\mathbb{R}
^{2d}$ be a relatively separated set of sampling points in $\mathbb{R}^{2d}$. Let $m: \mathcal{X}\rightarrow \mathcal{J}$ be a mapping assigning a window $g^{m(z_i)}\in M $ to each sampling point. 
Then the set 
\begin{equation}\label{BSgen}
	\{\pi (z_i)g^{m(z_i)}, i\in\mathcal{I}\} 
\end{equation}
is a Bessel Sequence for $\mathbf{L}^2 (\mathbb{R}^d )$. 
\end{Th}
 \underline{Proof:}\\ We have to estimate the series
 \[\sum_{z_i\in \mathcal{X}}|\langle f ,
 \pi (z_i ) g^{m(z_i )}\rangle |^2 =
  \sum_{z_i\in \mathcal{X}}|\mathcal{V}_{g^{m(z_i )}} f (z_i)|^2
   \leq \sum_{ z_i \in \mathcal{X}} |F(z_i)|^2
   ,\]
with  $F$ given by
(\ref{point-max1}).
 Now, as shown in the previous lemma, $F$ is in $\Wsp(C^0,\ell^2 )$. Hence, 
since  $\mathcal{X}$  is a relatively separated set,
Lemma~\ref{Le:supest} can be applied to obtain the following
estimate  for
the Bessel  bound of \eqref{BSgen}:
\[\sum_{z_i\in \mathcal{X}}|\langle f ,
 \pi (z_i ) g^{m(z_i )}\rangle |^2\leq C^2 R(1+1/\gamma)^{2d}\|f\|^2_2
\]
\hfill$\square$

\begin{Lem}{\label{union-relsep} }
The union of points $\{\lambda =
(x,\xi ), \lambda\in\mathcal{X}^r\}$ in the  discrete sets $\mathcal{X}^r$  as
chosen in Theorem~\ref{mainth} is relatively separated.
\end{Lem}
\underline{Proof:}\\Each of the lattices $\Lambda^r$ determining the
 Gabor frames used in the construction of
$\Gmgr$ is of course separated, even uniformly with respect to
$r$.  The admissibility condition for  $(\Omega _r
)_{r\in\mathcal{I}}$ allows only finite overlap between any pair of members
in the covering, or equivalently that the family of  balls of radius $R_r$
centered at $x_r$ forms a covering of (uniformly) bounded height.
By assumption, increasing each of the balls $B_{R_r}(x_r)$ by the
finite radius $\delta >0$, only the height of the covering will be
increased, but not the property of (uniformly) finite height. In
other words, the family of enlarged balls $B_{R_r + \delta}(x_r)$ is
still an admissible covering of ${\mathbb{R}^d}$ and the union of
the $\mathcal{X}$ is a relatively separated set.  \hfill$\square$\\[5mm]

We conclude the proof of Theorem~\ref{mainth}  by choosing $ M = \{g^j, j \in \mathcal{J}\}$ in the following Corollary. 
\begin{Cor}\label{cor:uppb}
An upper frame bound for $\Gmgr$ as defined in  (\ref{mixsyst1})
is given by  $C^2 n(\delta )(1+1/\gamma)^{2d}$.
\end{Cor}
Note that $n(\delta )$ denotes the height of the covering, which depends on $\delta$.
\section{Reduced multi-window Gabor frames:  windows with compact support or bandwidth}
In this section, we show that in a specific situation, which is, however, of practical relevance, quilted Gabor frames may be constructed. In the present model, we only change the resolution in  time (or frequency). This means, that the time-frequency domain is partitioned into stripes rather than patches. Under the additional assumption that the analysis window has compact support (or bandwidth), we easily obtain a lower frame bound for the quilted system. \\

Assume that we are given  Gabor frames
$\mathcal{G}(g^j,\Lambda^j )$ for $\Lz$, $j \in \mathcal{J}$,  where all the
windows $g^j$ have compact support and $\|g^j\|_{S_0}\leq C_g<\infty \,\forall j$. We now want to use each system for a certain time, i.e., in a restricted stripe in the time-frequency domain.  The stripes are defined 
   by means of
a partition of unity:  we assume that
$f = \sum_{r\in\mathcal{I}} \psi_r f $  with $\psi_r\leq 1$ and that the $\psi_r$ 
have compact support  in $[a_r, b_r]$. By means of a mapping $m:\mathcal{I}\rightarrow \mathcal{J}$, we assign one particular frame  to each of these stripes.

Now, subfamilies of the given Gabor frames may be  constructed as
follows. Assume, for simplicity, that $m(0) = 0$ and consider the task to represent $\psi_0 f$ in terms of the given
Gabor frame $\mathcal{G}(g^0,\Lambda^0 )$:
\begin{align*}
\psi_0 f  =& \psi_0(\sum_{\lambda \in \Lambda^0} \langle f, \pi
(\lambda )g^0\rangle \pi (\lambda )g^0)\\
=&\sum_{\lambda \in \Lambda^0} \langle f, \pi (\lambda )g^0\rangle
\psi_0 \pi (\lambda )g^0.
\end{align*}
Now,  there exist $n^u_0$ and $n^l_0$ such that for $\lambda = (n\alpha_0, m\beta_0 )$ with $n\notin [n^l_0, n^u_0]$, we find that $\psi_0  \pi (\lambda )
g^0\equiv 0$, hence
\[\psi_0 f  = \sum_{\lambda \in \mathcal{X}_0} \langle f, \pi (\lambda )g^0\rangle
\psi_0\pi (\lambda ) g^0 ,\]where $\mathcal{X}_0 = [n^l_0\cdot \alpha_0, n^u_0\cdot \alpha_0]\times\beta_0\mathbb{Z}$ is the subset of
$\Lambda^0$ corresponding to the nonzero contributions.

Analogously  subsets $\mathcal{X}_r\subset \Lambda^{m(r)}$ are
chosen for all $r$, and we obtain:
\begin{equation}\label{Eq:RedMulWin}
f = \sum_{r\in\mathcal{I}} \psi_r f = \sum_r \sum_{\lambda \in \mathcal{X}^r} \langle f, \pi (\lambda )g^{m(r)}\rangle
\psi_r\pi (\lambda ) g^{m(r)} .
\end{equation} With this construction, we  state the following proposition.
\begin{Pro}\label{Pr:QuiCp} 
For a family of tight Gabor frames $\mathcal{G}(g^j,\Lambda^j )$,  $j \in \mathcal{J}$,  for $\Lz$ let 
$\sup_{j\in\mathcal{I}}\|g^j\|_{S_0} = C_g<\infty	$
 and $C_{\Lambda^j} = (\frac{1}{\alpha_j}+1)^{d}(\frac{1}{\beta_j}+1)^{d} \leq C_{\Lambda}<\infty$ for all $j\in\mathcal{J}$. Let a partition of unity
$(\psi_r )_{r\in\mathcal{I}}$  of compactly supported $\psi_r$ with height $n_0$ be given
and let  a mapping $m:\mathcal{I}\rightarrow \mathcal{J}$ assign a frame  $\mathcal{G}(g^{m(r)} ,\Lambda^{m(r)} )$ to each $r\in\mathcal{I}$.
Assume that index sets 
$\mathcal{X}^r   = [n^l_r\cdot \alpha_{m(r)} , n^u_r\cdot \alpha_{m(r)}]\times\beta_{m(r)}\mathbb{Z}$ are chosen such that 
for all $r\in\mathcal{I}$ and $\lambda = (n\alpha_{m(r)}, m\beta_{m(r)} )$ with $n\notin [n^l_{m(r)}, n^u_{m(r)}]$, we have that $\psi_r  \pi (\lambda )
g^{m(r)}\equiv 0$.\\
Then, 
the union of the subfamilies $\bigcup_{r\in\mathcal{I}}(g^{m(r)},\mathcal{X}_r )$ is a
frame for $\Lz$ with a lower frame bound given by $1/(n_0 C_{\Lambda}C^2_g)$.
\end{Pro}
\underline{Proof:}\\
First note that
\begin{equation}\label{one}
\|\psi_r h\|_2^2\leq \|h\|_2^2 \ \mbox{ for all }\ h\in\Lz
\end{equation}
 Now set $h_r = \sum_{\lambda \in \mathcal{X}^r} \langle f, \pi (\lambda )g^{m(r)}\rangle
\pi (\lambda ) g^{m(r)}$ and thus, with \eqref{one}:
\begin{align*}
  \| f\|_2^2\leq &n_0 \sum_r\|\psi_r f\|_2^2 \\
 =& n_0 \sum_r\| \sum_{\lambda \in \mathcal{X}^r} \langle f, \pi (\lambda )g^{m(r)}\rangle
\psi_r \pi (\lambda )g^{m(r)}\|_2^2\\
=&n_0 \sum_r\| h_r \psi_r\|_2^2\leq  n_0 \sum_r\|\sum_{\lambda \in \mathcal{X}^r} \langle f, \pi
(\lambda )g^{m(r)}\rangle \pi (\lambda ) g^{m(r)}\|_2^2\\
\leq & n_0 \sum_r (\frac{1}{\alpha_j}+1)^{d}(\frac{1}{\beta_j}+1)^{d} \|g^{m(r)}\|^2_{S_0}
\sum_{\lambda \in
\mathcal{X}^r} |\langle f, \pi (\lambda )g^{m(r)}\rangle|^2
\\
\leq& n_0 C_{\Lambda}C^2_g\sum_r\sum_{\lambda \in
\mathcal{X}^r} |\langle f, \pi (\lambda )g^{m(r)}\rangle|^2
\end{align*}
The last inequality is due to the boundedness of the
frame-synthesis operator $T^\ast_{g^{m(r)}}: \ell^2 (\Lambda^j ) \rightarrow \Lz$, whenever the window $g^{m(r)}$
is in $S_0 (\mathbb{R}^d )$, see~\cite[Proposition~6.2.2]{gr01}. This proves the existence of a lower frame bound as stated. The existence of an upper frame bound can be seen directly from the construction of the subfamilies, and is furthermore  covered by the general 
case proved in Section~\ref{Se:BesselCond}.
\hfill$\square$
\begin{rem}
\begin{enumerate}
\item An analogous statement holds for general, not necessarily tight frames, for,  if $\gamma^j$ are the dual windows for each $g^j$, then 
$\|S_j^{-1} g^j\|\leq \frac{1}{A_j}\|g^j\||$.
\item The same construction may be realized in the Fourier transform domain by applying a partition of unity to $\hat{f}$. This corresponds to the usage of different Gabor frames in different stripes of the frequency domain and hence resembles a non-orthogonal filter bank. As a particular example, a constant-Q transform may be realized~\cite{BR91}.
\item Note that a similar yet more restrictive construction, corresponding to the classical  ``painless non-orthogonal expansions" was suggested in~\cite{badoja09}. Very recently, another related and highly interesting construction has been suggested in~\cite{prruwo09}.
\end{enumerate}
\end{rem}
%%%%%%%%%%%%%%%%%%%%%%%%%%%%%%%%%%%%%%%%%%%%%%%%%%%%%%%%%%%%%%%%%%%%%%%%%%%%%%%%%%%%%%%%%%%%%%%%%%%%%%%%%%%%%%%%%%%%%%%%%%%
\section{Replacing a finite number of  frame elements }

We now consider the task of replacing a  finite number of atoms from  a given (Gabor) frame  by  a finite number of atoms from a different (Gabor) frame. The following theorem gives a condition valid for general frames, which will then be applied to Gabor frames.  Recall, that $T_g$ and $T_g^\ast$ denote the analysis and synthesis mapping, respectively, for   given $g$ and $\Lambda$. In this section, we use the notation $T_{g,\mathcal{I}}$ and $T^\ast_{g,\mathcal{I}}$ for the respective mappings corresponding to subsets of the given lattices. For example, let $\mathcal{F}_1\subset\Lambda$ be a finite subset of a lattice $\Lambda$, then $T_{g,\mathcal{F}_1}(f) = \langle f,\pi (\lambda ) g\rangle$, for $\lambda\in\mathcal{F}_1$. The theorem makes use of a linear mapping $L$ describing the replacement procedure in the coefficient domain. As long as  elements from frame $\mathcal{G}_1$ may be replaced by elements from $\mathcal{G}_2$ in a controlled manner, i.e. without loosing energy, a quilted frame can be obtained.   
\begin{Th}\label{Th:RepFin}Assume that two frames $\mathcal{G}_1= \{g_i, i\in\mathcal{I}\} $ and $\mathcal{G}_2=\{h_j, j\in\mathcal{J}\} $ for $\mathbf{L}^2  (\mathbb{R}^d)$ are given and a finite number of elements  $\{g_i, i\in\mathcal{F}_1\subset\mathcal{I}\}$ of $\mathcal{G}_1$ are to be replaced by a finite number of elements  $\{h_j, j\in\mathcal{F}_2\subset\mathcal{J}\}$ of $\mathcal{G}_2$. \\
Let $A_1$ be the lower frame-bound of $\mathcal{G}_1$.
If a bounded linear mapping $L:\ell^2 (\mathcal{F}_1 )\mapsto \ell^2 (\mathcal{F}_2 )$ can be found such that 
\begin{equation}\label{Eq:RepCond}
\| T^\ast_{g,\mathcal{F}_1} - T^\ast_{h,\mathcal{F}_2}L\|_2^2 = C<\frac{A_1}{2},
\end{equation} then the set
\begin{equation}\label{Eq:rep1}\{g_i, i\in\mathcal{I}\backslash\mathcal{F}_1\}\cup\{h_j, j\in\mathcal{F}_2\}
\end{equation}
 is a frame for $\mathbf{L}^2 (\mathbb{R}^d )$ with a lower frame bound given by $(A_1-2C)/\max(1,2\|L^\ast \|_2^2)$.
\end{Th}
\underline{Proof:}\\
First note that 
\[A_1\|f\|_2^2\leq\sum_{i\in\mathcal{I}\backslash\mathcal{F}_1}|\langle f,g_i\rangle |^2+\sum_{i\in\mathcal{F}_1}|\langle f,g_i\rangle |^2 = \|T_{g,\mathcal{I}\backslash\mathcal{F}_1}\|_2^2+\|T_{g,\mathcal{F}_1}\|_2^2 .\]
Now, we have
\begin{align*}
	\|T_{g,\mathcal{F}_1}\|_2^2 &\leq (\|(T_{g,\mathcal{F}_1}-L^\ast T_{h,\mathcal{F}_2})f\|+\|L^\ast T_{h,\mathcal{F}_2}  f\|)^2\\
	&\leq 2\cdot \|(T_{g,\mathcal{F}_1}-L^\ast T_{h,\mathcal{F}_2})f\|_2^2+2\cdot \|L^\ast T_{h,\mathcal{F}_2}\|_2^2\\
	&\leq 2\cdot C\cdot \|f\|_2^2+2\cdot \|L^\ast T_{h,\mathcal{F}_2}f\|_2^2,
\end{align*}
hence
\begin{align*}
(A_1-2C)\|f\|_2^2&\leq\|T_{\mathcal{I}\backslash\mathcal{F}_1}f\|_2^2+2\cdot \|L^\ast \|_2^2\|T_{h,\mathcal{F}_2}f\|_2^2 \\
&\leq \max(1,2\|L^\ast \|_2^2)\cdot
 (\sum_{i\in\mathcal{I}\backslash\mathcal{F}_1}|\langle f,g_i\rangle |^2 +\sum_{j\in\mathcal{F}_2}|\langle f,h_j\rangle |^2)
 \end{align*} and hence $(A_1-2C)/\max(1,2\|L^\ast \|_2^2)$ is a lower frame bound for the system given in \eqref{Eq:rep1}. The existence of an upper frame bound is trivial. \hfill$\square$\\[5mm]
Note that the above theorem only states the existence of a lower frame bound under the given conditions, while this frame bound will usually not be optimal.\\
We now turn to  the special case of Gabor frames.

\begin{Cor}\label{Cor:GabRep}
Assume that  tight Gabor frames $\mathcal{G}(g,\Lambda^1)$ and $\mathcal{G}(h,\Lambda^2)$ with
$g,h\in S_0(\mathbb{R}^d )$ are given. Assume further that in a compact region
$\Omega\subset\mathbb{R}^{2d}$ the time-frequency shifted atoms $\pi
(\lambda ) g$, $\lambda\in\mathcal{F}_1 = \Omega\cap\Lambda^1$ are
to be replaced by a finite set of time-frequency shifted atoms $\pi (\mu
) h$, $\mu\in\mathcal{F}_2\subset \Lambda^2$. 
\begin{itemize}
\item[(a)] If $\Lambda^1 = \Lambda^2$, then 
\begin{equation}
\mathcal{G} =\{\pi(\lambda ) g:
\lambda\in\Lambda\setminus\mathcal{F}_1\}\cup \{\pi(\mu) h: \mu\in\mathcal{F}_2\}\end{equation}is a (quilted Gabor) frame for $\mathbf{L}^2(\mathbb{R}^d)$, whenever 
\[\|h-g\|_{S_0}^2\ = C<\frac{A_1}{2C_\Lambda},\]where $C_\Lambda$ is a constant only depending on  the lattice $\Lambda = \Lambda^1 = \Lambda^2$.
\item[(b)]For general $\Lambda^2$, a compact set $\Omega^{\ast}$ in  $\mathbb{R}^{2d}$ can be chosen such that for 
$\mathcal{F}_2 =\Omega^{\ast}\cap\Lambda^2$, the union 
\begin{equation}\label{Eq:Rep1}
\mathcal{G} =\{\pi(\lambda ) g:
\lambda\in\Lambda^1\setminus\mathcal{F}_1\}\cup \{\pi(\mu) h: \mu\in\mathcal{F}_2\}\end{equation}
is a (quilted Gabor) frame.
\end{itemize} 

\end{Cor}
\underline{Proof:}\\
Statement (a) can easily be seen by choosing $L = \Id$ in  Theorem~\ref{Th:RepFin}:
\begin{align*}
	\|T_{\mathcal{F}_1}^\ast-T_{\mathcal{F}_1}^\ast \|^2_{\ell^2\mapsto \mathbf{L}^2(\mathbb{R}^d) }
	&=\sup_{\|c\|_2 = 1} \|\sum_{\lambda\in\mathcal{F}_1}c_\lambda (\pi (\lambda) g - \pi (\lambda) h)\|^2_2\\
	&\leq C_\Lambda\cdot \|h-g\|_{S_0}^2,  
\end{align*}where the last inequality follows from boundedness of the synthesis operator for windows in $S_0$,~\cite[Theorem~12.2.4]{gr01}.\\
To show (b), we first introduce the mapping $L:\ell^2 (\mathcal{F}_1 )\mapsto \ell^2 (\mathcal{F}_2 )$ as follows. For a finite sequence $\mathbf{c} = (c_\lambda)_{\lambda\in\mathcal{F}_1}$ we define $L(\mathbf{c} )(\mu )  = (\langle\sum_{\lambda\in\mathcal{F}_1}c_\lambda \pi (\lambda) g , \pi (\mu) h\rangle)_{\mu\in\mathcal{F}_2}$, for which 
\begin{align*}
\|L\|_{\ell^2 (\mathcal{F}_1 )\mapsto \ell^2 (\mathcal{F}_2 )}^2 &=\sup_{\|c\|_2 = 1}\sum_{\mu\in\mathcal{F}_2} |\sum_{\lambda\in\mathcal{F}_1}
c_\lambda \langle \pi (\lambda) g ,\pi (\mu) h\rangle|^2\\
&= \|T_{h,\mathcal{F}_2}T^\ast_{g,\mathcal{F}_1}\|_2^2\leq C_{\Lambda^1}C_{\Lambda^2}\cdot\|g\|_{S_0}^2\|h\|_{S_0}^2 
\end{align*}
due to the boundedness of synthesis and analysis operator, $T_{h,\mathcal{F}_2}$ and $T^\ast_{g,\mathcal{F}_1}$, respectively. We then have:
\begin{align*}
	\|T_{\mathcal{F}_1}^\ast-T_{\mathcal{F}_2}^\ast L\|_{\ell^2\mapsto \mathbf{L}^2(\mathbb{R}^d) } &=\sup_{\|c\|_2 = 1} \|\sum_{\lambda\in\mathcal{F}_1}c_\lambda \pi (\lambda) g -\sum_{\mu\in\mathcal{F}_2} \sum_{\lambda\in\mathcal{F}_1}\langle c_\lambda \pi (\lambda) g,\pi (\mu) h\rangle \pi (\mu) h\|_2\\
%	&=\sup_{\|c\|_2 = 1} \|\sum_{\mu\in\Lambda^2}\sum_{\lambda\in\mathcal{F}_1}c_\lambda \langle \pi (\lambda) g,\pi (\mu) h\rangle \pi (\mu) h -\sum_{\mu\in\mathcal{F}_2} \sum_{\lambda\in\mathcal{F}_1}c_\lambda \langle \pi (\lambda) g,\pi (\mu) h\rangle \pi (\mu) h\|_2\\
	&=\sup_{\|c\|_2 = 1} \|\sum_{\mu\in\mathcal{F}_2^c}\sum_{\lambda\in\mathcal{F}_1}c_\lambda \langle \pi (\lambda) g,\pi (\mu) h\rangle \pi (\mu) h \|_2\\
	&\leq \|h\|_2\sup_{\|c\|_2 = 1} \sum_{\mu\in\mathcal{F}_2^c}\sum_{\lambda\in\mathcal{F}_1}|c_\lambda |\cdot |\mathcal{V}_g h (\mu-\lambda ) |\\
	&\leq \|h\|_2\sqrt{|\mathcal{F}_1| } \sum_{\mu\in\mathcal{F}_2^c} \max_{\lambda\in\mathcal{F}_1}|\mathcal{V}_g h (\mu-\lambda ) |
\end{align*}
where $\mathcal{F}_2^c = \Lambda^2\backslash\mathcal{F}_2$.
Now an appropriate set  $\mathcal{F}_2$ may be constructed as follows: \\
Let $G = |\mathcal{V}_{g}h|\in \Wsp(C^0,\ell^1 )$. We can choose a compact set $\Omega^\ast\subset\mathbb{R}^{2d}$, such that
$$\|\max_{\lambda\in\mathcal{F}_1}(T_\lambda G- (T_\lambda G)\cdot\chi_{\Omega^\ast})\|_{W} < \tilde{\varepsilon}=\frac{\sqrt{A_1}}{\sqrt{2|\mathcal{F}_1|}\|h\|_2 C_{\Lambda^2}} ,$$ where $\chi_{\Omega^\ast}$ is the indicator function of the set $\Omega$. 
Now set 
 \[\mathcal{F}_2 = \Lambda^2\cap\Omega^\ast ,\]then
 \begin{align}\notag
\sum_{\mu\in\mathcal{F}_2^c}\max_{\lambda \in\mathcal{F}_1}
T_{\lambda}G (\mu) &= \| (\max_{\lambda \in\mathcal{F}_1}T_{\lambda}G) \big|_{\mathcal{F}_2^c}\|_{\ell^1}\\\notag
&=\|  \max_{\lambda\in\mathcal{F}_1}(T_\lambda G- T_\lambda G\cdot\chi_{\Omega^\ast})  \big|_{\Lambda^2}\|_{\ell^1}\\\label{eqwamal}
&\leq C_{\Lambda^2}\cdot \|(T_{\lambda}G-T_{\lambda}G\cdot\chi_{\Omega^\ast})\|_{\Wsp(C^0,\ell^1)}\\
&\leq  C_{\Lambda^2}\cdot\tilde{\varepsilon}= C_{\Lambda^2}\cdot\frac{\sqrt{A_1}}{\sqrt{2|\mathcal{F}_1|}\|h\|_2 C_{\Lambda^2}}, \notag
\end{align}
where we used Lemma~\ref{amalest1} in \eqref{eqwamal}.
Hence, $\mathcal{F}_2$ can be chosen, such that 
\[\|T_{\mathcal{F}_1}^\ast-T_{\mathcal{F}_2}^\ast L\|^2_{\ell^2\mapsto \mathbf{L}^2(\mathbb{R}^d) } <A_1/2.\] Of course, the existence of an upper frame bound is trivial
and the resulting system \eqref{Eq:Rep1} is a quilted Gabor frame according to Theorem~\ref{Th:RepFin}.\hfill$\square$\\
 %with lower frame bound $\frac{A_1-2\sqrt{|\mathcal{F}_1|}C_{\Lambda^2}\varepsilon}{\max (1,2\cdot B_1\cdot\|h\|_2^2 )}$. 
 
%Hence, we also have $\|T_\lambda(G-G\chi_{\Omega})\|_{W} < \varepsilon$. 
 %We have:\\ If $g_1$, $g_2$
%are in $S_0$, then the short-time Fourier transform
%$|\mathcal{V}_{g_2}g_1|$ is in $\Wsp(C_0,\ell_1)$,\\
%and:\\
%If $G\in \Wsp(C_0,\ell_1)$ and the set $\mathcal{X}$ is separated,
%then for all $\varepsilon> 0$, there exists a compact set
%$\Omega$, such that
%\[\sum_{x_i\in\mathcal{X}\cap\Omega^c} G(x_i ) < \varepsilon\]
%Together:\\
%Find a compact set

\begin{rem}
(1) The size of the region in which atoms from $\mathcal{G}_1$ have to be replaced, influences the choice of $\Omega^\ast$.  In particular, $\tilde{\varepsilon}$ is  reciprocally related to $\sqrt{|\mathcal{F}_1|}$, i.e., $\Omega^\ast$ grows in dependence on the perimeter rather than the area of $\Omega$.\\
(2) For statement (a) in Corollary~\ref{Cor:GabRep}, if  $g,h\in \mathbf{L}^2(\mathbb{R}^d)$, then the frame property follows whenever \[\|h-g\|_2^2\ = C<A_1/(2|\mathcal{F}_1|).\]
Note that  the size of $\Omega$ determines the necessary similarity of the windows, whereas, for windows in $S_0(\mathbb{R}^d)$, the good localization implied by $S_0$-membership allows for a stronger statement. \\
(3) The construction in Corollary~\ref{Cor:GabRep} implies  explicit dependence on time and frequency of the resulting quilted frame. Similarly, Gabor atoms in finitely many compact areas can be replaced by different Gabor systems. Details on this procedure will be reported elsewhere. From an application point of view, this process corresponds to finding optimal representations for local signal components, e.g. in the sense of sparsity.  \\
(4) In Example~\ref{Ex00}, the linear mapping $L:\Lambda^0\rightarrow \bigcup_r\mathcal{X}^r$ can be chosen as follows: 
\[L(c)(\lambda) =\begin{cases}c_\lambda , & \mbox{for } \;\;\lambda \in \Lambda^0 \\ 0, & \mbox{for} \;\;\lambda \in \bigcup_r\mathcal{X}^r\backslash\Lambda^0 \end{cases}.\]
Then
\[\|T_{\Lambda^0,g}^\ast-T_{\bigcup_r\mathcal{X}^r,g}^\ast L  \|^2
	=\sup_{\|c\|_2 = 1} \|\sum_{\lambda\in\Lambda^0}c_\lambda \pi (\lambda) g - \sum_{\lambda\in\bigcup_r\mathcal{X}^r}L(c)(\lambda )\pi (\lambda) g\|^2_2 = 0.\]
As underlined by Example~\ref{Ex00} and Corollary~\ref{Cor:GabRep}(a), in constructing quilted Gabor frames, technically difficult situations mainly arise if the lattices \emph{and} the windows vary.
\end{rem}
\section{Reconstruction and Simulations}\label{Se:Mult-gab-rec}
Since the frame property has been proved for systems as described in Proposition~\ref{Pr:QuiCp} and Corollary~\ref{Cor:GabRep}, reconstruction can always be performed by means of a dual frame. However, since quilted Gabor frames possess a strong local structure, alternative and numerically cheaper methods may be preferable as long as sufficient precision in the reconstruction may be guaranteed. The next two sections present numerical results for various approaches to reconstruction for which the calculation of an exact dual frame is not necessary.   
\subsection{Reduced multi-window Gabor frames}
We first consider reduced multi-window Gabor frames. Here, equation~\eqref{Eq:RedMulWin} yields an immediate reconstruction formula by means of projections onto the members of the partition of unity. However, we are more interested in the generic reconstruction by means of dual frames. We may compare the dual frame corresponding to  the quilted Gabor frame    $\bigcup_{j\in\mathcal{I}}(g^j,\mathcal{X}_j )$ with the quilted system $\bigcup_{j\in\mathcal{I}}(\gamma_j,\mathcal{X}_j )$  resulting from using the dual windows $\gamma_j$ of the original frames $\mathcal{G}_j$. Alternatively, we may start with tight Gabor frames.\\  While this approach  does not result in perfect reconstruction in one step, we apply the frame algorithm (see~\cite[Chapter~5]{gr01}) to obtain near-perfect  reconstruction in  a few iteration steps. In this context, the condition number of the operators involved plays an important  role  and depends on the amount of overlap that we introduce in the design of the system.  In the following example, it turns out that, while no essential overlap is necessary to obtain a frame in the finite discrete case, the overlap as requested in the proof of Proposition~\ref{Pr:QuiCp} leads to faster convergence of the frame algorithm.
\begin{ex}\label{Ex1}
\upshape We consider two tight Gabor frames for $\mathbb{C}^L$ and $L = 144$ and two cases of different redundancy. First, redundancy is $4.5$, corresponding to the lattices $\Lambda^1$ with $a = 4,b = 8$ and $\Lambda^2$ with $a = 8,b = 4$. Second, we consider two frames with redundancy $1.125$, corresponding to the lattices $\Lambda^1$ with $a = 8, b = 16$ and $\Lambda_2$ with $a = 16, b = 8$. The corresponding tight windows are shown in Figure~\ref{Fi:TW}. 
\begin{figure}[tb]
\centerline{\includegraphics[width = 10cm,height=5cm]{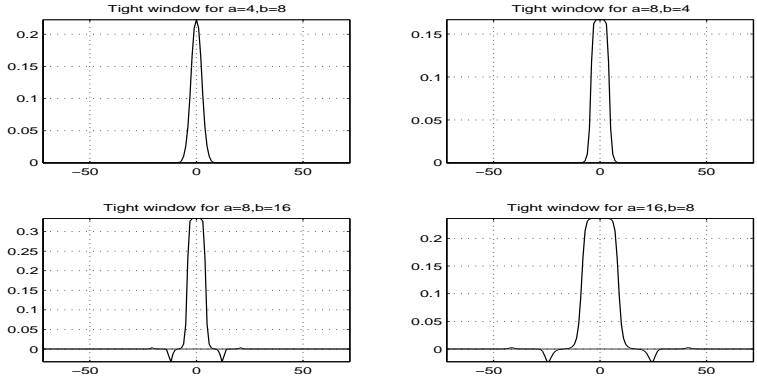}}
\caption{{\it Tight windows corresponding to four different Gabor frames }}\label{Fi:TW}
\end{figure}
We next generate a quilted Gabor system without overlap and a corresponding Gabor system with overlap, for both cases of redundancy. Note that in each lattice point as depicted in Figure~\ref{Fi:Overlaps}, the tight window of the original Gabor frame is used. 
\begin{figure}[tb]
\centerline{\includegraphics[scale = .6]{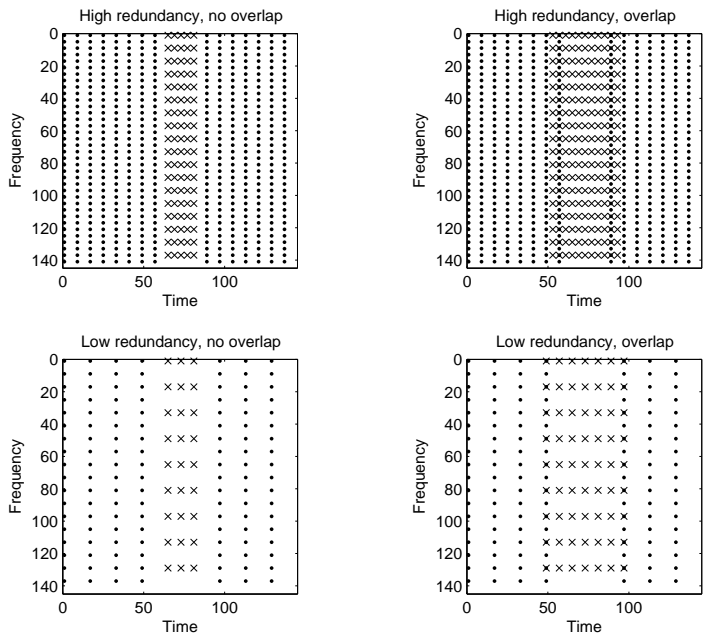}}
\caption{{\it Various lattices with different amounts of overlap }}\label{Fi:Overlaps}
\end{figure}
We now look at the condition numbers of the resulting quilted systems,  listed in Table~1.    
 \begin{table}\label{Tab1}
 \begin{center}{\bf Table~1: Condition number of four quilted Gabor frames}
\begin{tabular}[t]{|c||c|c|}
    \hline 
    Redundancy & No overlap &  Overlap \\\hline
    $4.5 $& $1.6$& $1.4$\\
    $1.125$& $6.4$& $1.5$\\\hline
  \end{tabular}
\end{center}\end{table}
It is obvious, that higher redundancy leads to more stability in the process of quilting frames. On the other hand, for the system with low redundancy, overlap becomes essential in order to obtain acceptable condition numbers. These observations are consistently confirmed by more extensive numerical experiments.\\
In a next step, we now compare the convergence of the (iterative) frame algorithm for the $4$ cases considered in this example. Table~2  gives the number of iterations necessary to attain the threshold of $10^{-8}$. Figure~\ref{Fi:ConvFi} then shows the rate of convergence for the three cases with acceptable conditions numbers.\\
 \begin{table}\label{Tab2}
 \begin{center}{\bf  Table~2: Number of iterations for convergence of frame algorithm}
\begin{tabular}[t]{|c||c|c|}
    \hline 
    Redundancy & No overlap &  Overlap \\\hline
    $4.5 $& $21$& $17$\\
    $1.125$& $322$& $18$\\\hline
  \end{tabular}
\end{center}\end{table}

\begin{figure}[tb]
\centerline{\includegraphics[scale = .5]{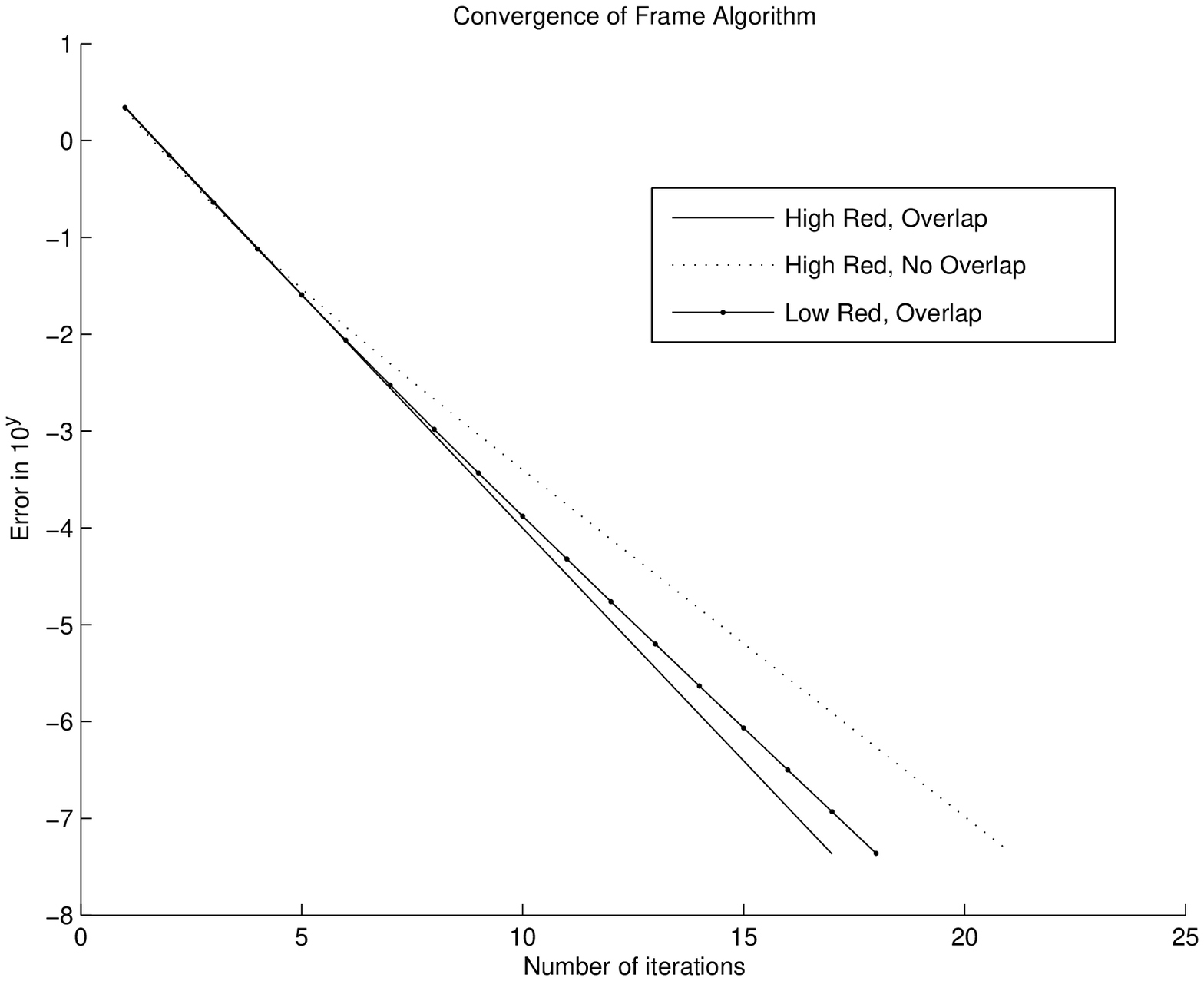 }}
\caption{{\it Convergence of frame algorithm} }\label{Fi:ConvFi}
\end{figure}
We finally discuss the following, ``preconditioned" version of reconstruction for the case of low redundancy with overlap. We wish to reconstruct a random signal $r\in\mathbb{C}^L$ from its quilted Gabor coefficient. As a first guess, instead of calculating the dual frame of the quilted frame, we simply use the quilted tight frame for reconstruction: Let $T_ {\mathcal{G}^t_q}$ denote the analysis operator corresponding to the quilted tight frame $\mathcal{G}^t_q=\bigcup_{j=1}^2(h_j,\mathcal{X}_j )$, where $h_j$ are the tight windows. We then obtain a reconstruction $rec$ by
\[rec = T_ {\mathcal{G}^t_q}^\ast \cdot T_ {\mathcal{G}^t_q}\cdot r.\]
 Obviously, the result is not accurate and in particular in the regions of transition between the two systems, errors occur. However, we can correct a considerable amount of the deviation from the identity by simply pre-multiplying the frame-operator  by the inverse of its diagonal. Hence, let
 $D$ be the diagonal matrix with $D(n,n)  = T_ {\mathcal{G}^t_q}^\ast \cdot T_ {\mathcal{G}^t_q}(n,n)$, $n = 1,\ldots , 144$. 
 \[rec = T_ {\mathcal{G}^t_q}^\ast \cdot T_ {\mathcal{G}^t_q}\cdot D^{-1}\cdot r.\]
  The respective results are shown in Figure~\ref{Fi:PreCondErr}. The relative error, defined by $\varepsilon = \frac{\|rec-r\|}{\|r\|}$ is then $0.2239$  for the uncorrected case and $0.032$ for the corrected version.
\begin{figure}[tb]
\centerline{\includegraphics[scale = .7]{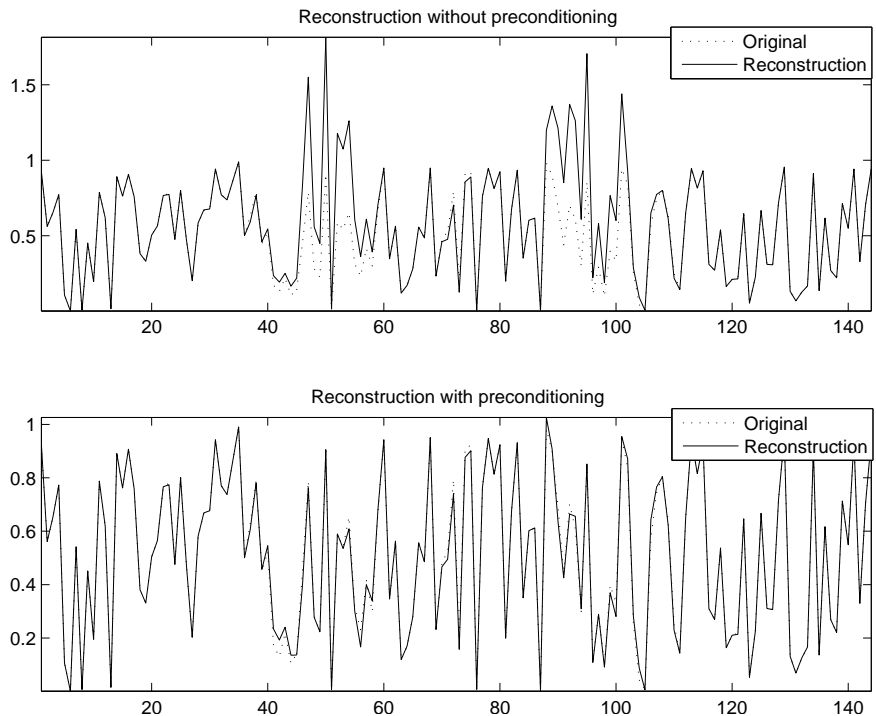}}
\caption{{\it Reconstructed random signal} }\label{Fi:PreCondErr}
\end{figure}
\end{ex}
\subsection{Replacing a finite number of elements}  
In our next example we consider a situation similar to the one discussed in Example~\ref{Ex1}, however, this time we wish to replace elements from $\mathcal{G}_1$ in a bounded, quadratic region of the time-frequency plane. 
\begin{ex}\label{Ex2} \upshape 
We consider the same Gabor frames as in Example~\ref{Ex1}, and look at the high redundancy systems first.  As before, we compare the condition number of the system obtained with overlap to the less redundant situation. The two situations are shown in Figure~\ref{Fi:Repl2Sys}. The quilted Gabor frame without overlap has condition number $1.5$, while allowing for some overlap, as shown in the second display of Figure~\ref{Fi:Repl2Sys} leads to condition number $1.4$. Accordingly, $18$  and $17$ iterations are necessary for convergence of the frame operator.\\
We now turn to the systems with low redundancy. Here we compare three amounts of overlap as shown in the upper displays of Figure~\ref{Fi:3Overlaps}. The condition numbers of the resulting systems and the convergence behavior of the corresponding frame algorithm are  shown  in the lower display of the same figure. Again, it becomes obvious that for low-redundancy systems, overlap is essential in order to obtain fast convergence in iterative reconstruction.  On the other hand, increasing overlap beyond a certain amount, does not dramatically improve the condition numbers. 
\begin{figure}[tb]
\centerline{\includegraphics[scale = .5]{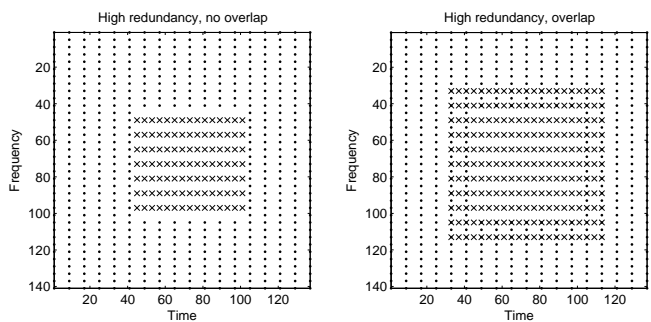}}
\caption{{\it Replacing atoms with and without overlap, high redundancy}}\label{Fi:Repl2Sys}
\end{figure}
\begin{figure}[tb]
\centerline{\includegraphics[scale = .7]{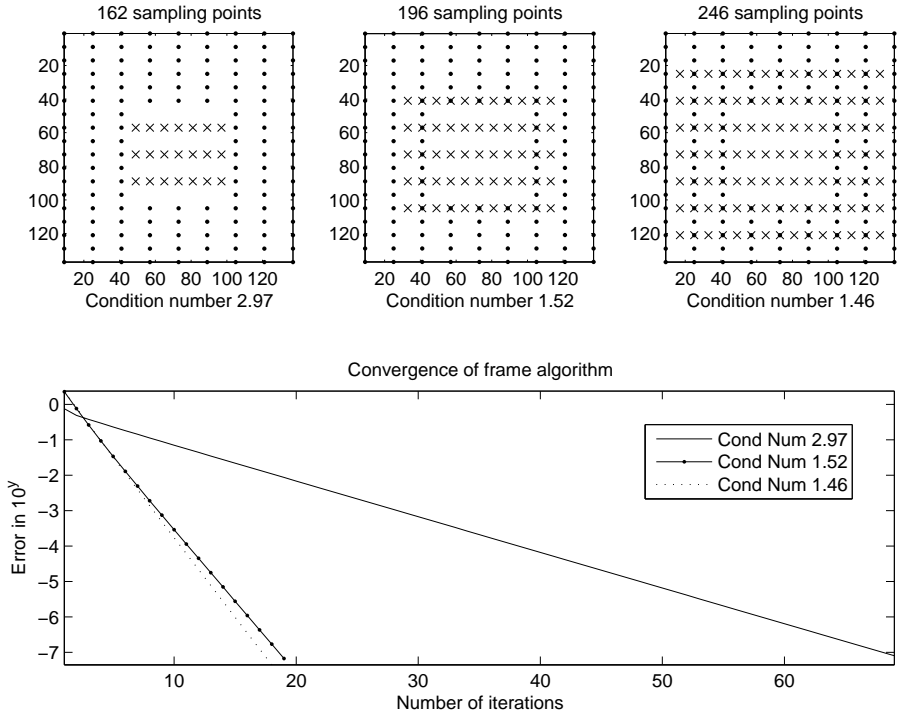}}
\caption{{\it Replacing atoms with and without overlap, low redundancy}}\label{Fi:3Overlaps}
\end{figure}

\end{ex}

\section{Summary and Outlook}\label{Se:Sum}
We have shown the existence of a lower frame bound for two particular instances of quilted Gabor frames. Furthermore, an upper frame bound has been constructed for the general setting. We showed how to reconstruct signals from the coefficients obtained with quilted Gabor frames and numerical simulations have been provided.\\
Future work will mainly include the construction of lower frame bounds for more general situations. In particular, Proposition~\ref{Pr:QuiCp} will be generalized to Gabor frames with non-compactly supported windows. Furthermore, numerical simulations suggest that atoms from a given Gabor frame may be replaced by atoms from a different frame in infinitely many compact regions of the time-frequency plane under certain conditions. 
On the other hand, for practical applications algorithms  applicable for long signals (number of sampling points $\gg 44100$) have to be developed.  The results of processing with quilted Gabor frames will be assessed on the basis of real-life data. Preconditioning similar to the procedure suggested in Example~\ref{Ex1} can be developed for the more complex situations of quilted frames, compare~\cite{bafehakr06-1}.
\section{Acknowledgments} 
The author wishes to thank Hans Feichtinger for the joint development of the notion of quilted frames as well as 	innumerable discussions on the topic, his proofreading and invaluable comments on the content of this article. She also wishes to thank Franz Luef for his comments on an earlier version of the paper and Patrick Wolfe for fruitful scientific exchange on the topic of adaptive frames from the point of view of applications.
%\bibliographystyle{abbrv}
%\bibliography{quiltbib,DT,general}
\def\cprime{$'$}

\end{document}